\input amstex
\documentstyle{amsppt}
\input bull-ppt
\keyedby{bull505e/paz}
\ratitle
\topmatter
\cvol{31}
\cvolyear{1994}
\cmonth{July}
\cyear{1994}
\cvolno{1}
\cpgs{64-67}
\title Alexander's and  Markov's theorems \\
in dimension four\endtitle
\subjclass Primary 57Q45, 57M25\endsubjclass 
\author Seiichi Kamada \endauthor
\address Department of Mathematics, Osaka City University,
Sumiyoshi, Osaka 558,
Japan\endaddress
\ml g54042\@jpnkudpc.bitnet\endml
\date May 10, 1993\enddate
\abstract Alexander's and Markov's theorems state that
any link type in $R^3$ is represented by a closed braid and
that such
representations are related by some elementary operations
called Markov moves.
We generalize the notion of a braid to that in 4-dimensional
space and
establish an analogue of these theorems.\endabstract
\endtopmatter


\document

\heading  1. Introduction \endheading
\par

    Any link type in $R^3$ is represented by a closed braid,
and two closed braids represent the same link type if and
only if they are related by braid isotopies, stabilizations,
and their inverse operations.   These
facts are well known as Alexander's theorem and Markov's
theorem
(cf. \cite{B; A; Mo}).
The latter is restated as follows: Two braids have ambient
isotopic closures in
$R^3$ if and only if  they are related by conjugations,
stabilizations, and their
inverse operations (Markov moves). These theorems enable
braid theory to play an
important role in (classical)  knot theory (for example,
\cite{J}).  We have a
natural analogue of them in 4-dimensional space.
\par

There seem to be many generalizations of the notion of
classical braids to
higher dimensions (\cite{D; G1; G2; BS; MS}, etc.).
We use, as a generalization of a classical braid,
the notion of a 2-dimensional braid due to O. Viro (cf.
\cite{K2}).  A similar notion was studied by L.~Rudolph
as a braided surface \cite{R1; R2}.
\dfn{Definition} A  2-{\it dimensional braid\/}
(of degree $m$)  is a compact oriented surface $B$
embedded in bidisk $D_1^2\times D_2^2$ such that (1) the
restriction of the second factor projection
 $p_2:D_1^2\times  D_2^2\to D_2^2$ to
$B$ is an oriented branched covering map of degree $m$
and (2) the boundary of $B$ is a trivial closed braid
$X_m\times\partial D_2^2 \subset
D_1^2\times\partial D_2^2 \subset  \partial(D_1^2\times
D_2^2)$, where $X_m$ is
fixed $m$ points in the interior of $D_1^2$. \enddfn

We do not require that the associated branched covering
map  $B\to D_2^2$ is {\it simple\/}, although it is
assumed in \cite{K2; R1; R2}.
\par

Two 2-dimensional braids are {\it equivalent\/} if they
are ambient isotopic by a fiber-preserving isotopy of
$D_1^2\times D_2^2$, where we
regard  $D_1^2\times
D_2^2$ as a $D_1^2$-bundle over $D_2^2$, keeping
$D_1^2\times\partial D_2^2$ fixed.  We do not distinguish
equivalent 2-dimensional braids.
Two 2-dimensional braids are {\it braid isotopic\/}
if we can deform one to the
other by an isotopy of $D_1^2\times D_2^2$ keeping the
condition of a
2-dimensional braid and $D_1^2\times\partial D_2^2$ fixed.
\par

Let $B$ and $B'$ be 2-dimensional braids of the same
degree $m$ in
$D_1^2\times D_2^2$ and in $D_1^2\times {D_2^2}'$.  Take a
boundary connected sum  $D_2^2\natural {D_2^2}'$ of
$D_2^2$ and ${D_2^2}'$ which
is also a 2-disk.  Then $B\cup B'$ in $D_1^2\times
(D_2^2\natural {D_2^2}')$ forms
a 2-dimensional
braid.  We call it a {\it braid sum\/} of $B$ and $B'$ and
denote it by $B\cdot
B'$, which is
uniquely determined up to equivalence.  The family of
equivalence classes of 2-dimensional braids of degree $m$,
together with this braid
sum operation, is a commutative semi-group. We call it
the 2-{\it dimensional braid semi-group\/}
of degree $m$.
\par

A 2-dimensional braid $B$ is naturally extended to a closed
oriented surface
$\widehat{B}$ embedded in
$D_1^2\times S^2= D_1^2\times (D_2^2\cup \overline{D_2^2})$
with
$\widehat{B}\cap D_1^2\times \overline{D_2^2}= X_m\times
\overline{D_2^2}$.
Identifying $D_1^2\times S^2$ with the tubular
neighborhood of the standard
2-sphere about the $t$-axis in $R^4$, we assume
$\widehat{B}$ to be  a surface
embedded in $R^4$ and call it the {\it closure\/} of
$B$ (about the $t$-axis)
in $R^4$ or a {\it closed {\rm2}-dimensional braid\/}.
\par

\proclaim{Theorem 1 \rm(Viro \cite{K1})}   Any closed
oriented surface embedded in $R^4$ is ambient isotopic to
a closed \RM2-dimensional braid in $R^4$.
\endproclaim
\par

  In \cite{K1} a closed 2-dimensional braid in $R^4$ is
treated as a sequence
of closed braids in classical dimension,  similar ideas to
which are found in
\cite{G; R1; K3}.
\par

Two operations, ``conjugations" and ``stabilizations",  are
defined for 2-dimensional braids as a natural 
analogue of
classical ones as
follows:
\par

Regard $D_1^2$ as $I_1^1\times I_2^1$ and $X_m$ as the set
$\{z_1,\cdots, z_m\}\subset D_1^2$ where $I_i^1$ $(i=1,2)$
is the interval $[-1,1]$ and  $z_i$
$(i=1,\cdots, m)$ is point $(0, 1-2^{-i})$.
\par

Let $b$ be a classical braid of degree $m$ in
$D_1^2\times I$ with $\partial b
=  b\cap D_1^2\times\partial I= X_m\times\partial I$.  The
product
$b\times S^1$  is a collection of annuli embedded in
$ D_1^2\times I\times S^1$. Identify $I\times S^1$ with the
collar
neighborhood  $N(\partial D_2^2)$ of $\partial D_2^2$ in
$D_2^2$ so that
 $\partial D_2^2\times\{1\}= \partial D_2^2$ and
$D_1^2\times I\times S^1$ with
$D_1^2\times N(\partial D_2^2)$.
\par

Up to equivalence, a 2-dimensional braid $B$ is assumed to
satisfy that
$B\cap  D_1^2\times N(\partial D_2^2)$ is the product
$X_m \times N(\partial D_2^2)$.
Replace it by $b\times S^1$ and we have a
2-dimensional braid.  It is called  a 2-dimensional
braid obtained from $B$ by the {\it conjugation\/}
associated with  $b$ or the
{\it conjugate\/} of $B$ by  $b$.
\par

For a nonnegative integer $a$, let $T_a$ be a homeomorphism
of $D_1^2= I_1^1\times I_2^1$ such that
$T_a(x, y_i)= (x, y_{i+a})$ for
each $x\in I_1^1$ and $y_i= 1-2^{-i}\in I_2^1$
$(i=0,1,2,\cdots)$, and let
$\widetilde{T_a}$ be the homeomorphism of
$D_1^2\times D_2^2$
which acts on the first factor by $T_a$ and on the second
factor by the identity.
For a 2-dimensional braid $B$ of degree $m$,
$\widetilde{T_a}(B)$
satisfies condition (1) of the definition of a 2-dimensional
braid and condition (2)
replaced $X_m$ by $T_a(X_m)= \{z_{a+1},\cdots, z_{a+m}\}$.
Up to equivalence, we
may assume that  $\widetilde{T_a}(B)$ is contained in
$I_1^1\times (y_a, y_{m+a+1})\times D_2^2
\subset I_1^1\times I_2^1\times D_2^2 =
D_1^2\times
D_2^2$. For non-negative integers $a$ and $b$,  let
$\iota_a^b(B)$ be the union
$\{z_1,\cdots, z_a\}\times D_2^2\cup \widetilde{T_a}(B)
\cup\{z_{m+a+1},\cdots, z_{m+a+b}\}\times D_2^2$,
which is a 2-dimensional braid of
degree $m+a+b$. $\iota_a^b$ is an injection from the 
2-dimensional braid
semi-group of degree $m$ into that of degree $m+a+b$.
\par

By \cite{K2} the equivalence class of a 2-dimensional braid
of degree $2$ is determined only by the number of branch
points.  Let $B^\ast$ be a unique, up to
equivalence, 2-dimensional braid of degree $2$ with two
branch points.  It is a
unique generator of the 2-dimensional braid semi-group of
degree $2$.
    A 2-dimensional braid $B'$ of degree $m+1$ is said to
be obtained from a
2-dimensional braid $B$ of degree $m$ by a {\it
stabilization\/} if $B'$ is\ 
equivalent to the braid sum $\iota_0^1(B)\cdot 
\iota_{1}^0(B^\ast)$.

\proclaim{Theorem 2}  Two \RM2-dimensional braids have 
ambient
isotopic
closures in $R^4$ if and only if  they are related by braid
isotopies,
conjugations, stabilizations, and their inverse operations.
\endproclaim
\par

\demo{Remark} A {\it simple\/} 2-dimensional braid is a
2-dimensional braid $B$
whose associated branched covering $B\to D_2^2$ is simple.
An advanced version of Theorem~1 is found in \cite{K1}:
Any closed oriented surface
embedded in $R^4$ is ambient isotopic to the closure
of a {\it simple\/} 2-dimensional
braid. Then we have a natural question which is open at
present: For two simple 2-dimensional braids with ambient
isotopic closures in $R^4$,
are they related
only through {\it simple\/} 2-dimensional braids by  braid
isotopies,
conjugations, stabilizations, and their inverse operations?
\enddemo
\par

We work in the piecewise linear category, and surfaces in 
4-space are assumed to
be embedded properly and locally flatly.
\par

\heading 2. Outline of the proof \endheading
\par

Our proof follows \cite{B}.
Let $\ell$ denote the $t$-axis of $R^4$ and
$\pi: R^4\to R^3$
the projection
along $\ell$.  We say that an oriented 2-simplex
$A$ in $R^4$  is {\it in general position with respect to\/}
$\ell$ if there is no 3-plane in
$R^4$ containing both $A$ and $\ell$.
Then $\pi(A)$ is an oriented 2-simplex
in $R^3$ which
forms, together with the origin of $R^3$, an oriented 
3-simplex.  
Using the orientation of $R^3$, we define the
{\it sign\/} of $A$ valued
in $\{\pm 1\}$.    Theorem~1 is shown as follows:
Let $F$ be a closed oriented surface
in $R^4$ and $K$ a {\it division\/} of $F$,
which is a certain kind of
tessellation by 2-simplices (not a triangulation but its
generalization).
    We may assume that each 2-simplex of $K$ is in general
position with
respect to $\ell$.   The number of 2-simplices of $K$ with
negative sign
is denoted by $h(K)$.  If $h(K)=0$, then $F$ is a closed 
2-dimensional braid
about $\ell$ in $R^4$.

\proclaim{Lemma 1 \rm(Sawtooth Lemma)} If there is a
\RM2-simplex $A$ of $K$ with
negative sign, then there is a sawtooth over $A$ of $F$.
\endproclaim
\par

A {\it sawtooth\/} is a family of 3-simplices in $R^4$
satisfying a nice
condition such that the surgery result of $K$ is a division
$K'$ of another
surface in $R^4$ with $h(K')= h(K)-1$.  By this lemma,
we can deform $F$ to a
closed 2-dimensional braid.
\par

As an analogue of \cite{B}, we can define some {\it
operations\/} which
transform a division of a surface in $R^4$ without changing
the isotopy type of
the  surface.  A {\it deformation chain\/} is a sequence of
divisions of
surfaces in $R^4$ connected by operations.
\par

\proclaim{Lemma 2}  Let $K, K'$ be a deformation chain with
$h(K)=h(K')\neq
0$.   Then there is a deformation chain  $K, K_1,\cdots, 
K_s,
K'$  for some
$s\geq 1$ such that  $h(K_i)<h(K)$ for $i$ $(i=1,\cdots, 
s)$.
\endproclaim
\par

\proclaim{Lemma 3}  Let $K, K', K''$ be a deformation chain
with
$h(K')>h(K)$ and $h(K')>h(K'')$. Then there is a deformation
chain  $K,
K_1,\cdots, K_s, K''$ for some $s\geq 1$ such that
$h(K_i)<h(K')$ for $i$ $(i=1,\cdots, s)$. \endproclaim
\par

By the argument in \cite{B}, we have Theorem~2 as a
consequence of Lemmas~2
and 3.\ The details will appear elsewhere.
\par


\Refs
\ra\key{MS}
\ref \key A \by J. W. Alexander
\paper A lemma on systems of knotted curves
\jour Proc. Nat. Acad. Sci. U.S.A.
\vol 9 \yr 1923 \pages 93--95 \endref

\ref \key B \by J. Birman
\book Braids, links and mapping class group
\publ Ann. of Math. Stud., vol. 82
\publaddr Princeton Univ. Press, Princeton, NJ
\yr 1974  \endref

\ref \key BS \by E. Brieskorn and K. Saito
\paper Artin Gruppen und Coxeter Gruppen
\jour Invent. Math.
\vol 17 \yr 1972 \pages 245--271 \endref

\ref \key D \by D. Dahm
\paper A generalization of braid theory
\paperinfo Ph.D. thesis
\jour Princeton Univ., Princeton, NJ, 1962  \endref

\ref \key G \by F. Gonz\'alez-Acu\~na
\paper A characterization of \RM2-knot groups
\jour preprint
\endref

\ref \key G1 \manyby D. L. Goldsmith
\paper The theory of motion groups
\jour Michigan Math. J.
\vol 28 \yr 1981 \pages 3--17 \endref

\ref \key G2 \bysame
\paper Motion of links in the \RM3-sphere
\jour Math. Scand.
\vol 50 \yr 1982 \pages 167--205 \endref

\ref \key J \by V. F. R. Jones
\paper A polynomial invariant for knots via von Neumann
algebras
\jour Bull. Amer. Math. Soc. (N.S.)
\vol 12 \yr 1985 \pages 103--111 \endref

\ref \key K1 \manyby S. Kamada
\paper A characterization of groups of closed orientable
surfaces in \RM4-space
\jour Topology \vol 33 \yr 1944 \pages 113--122
 \endref

\ref \key K2 \bysame
\paper Surfaces in $R^4$ of braid index three are ribbon
\jour J. Knot Theory Ramifications
\vol 1 \yr 1992 \pages 137--160 \endref

\ref \key K3 \bysame
\paper \RM2-dimensional braids and chart descriptions
\jour Topics in Knot Theory, Proceedings of the NATO
ASI on Topics in Knot Theory, Turkey, 1992
(M. E. Bozh\"uy\"uk, ed.), pp. 277--287\endref

\ref \key MS \by Yu. I. Manin and V. V. Schechtman
\book Arrangements of hyperplanes, higher braid groups and
higher Bruhat orders
\publ Adv. Stud. Pure Math., vol. 17
\publaddr Academic Press, Boston, MA
\yr 1986 \pages 289--308
\endref

\ref \key Mo \by H. R. Morton
\paper Threading knot diagrams
\jour Math. Proc. Cambridge Philos. Soc.
\vol 99 \yr 1986 \pages 247--260 \endref

\ref \key R1 \manyby L. Rudolph
\paper Braided surfaces and Seifert ribbons for closed 
braids
\jour Comment. Math. Helv.
\vol 58 \yr 1983 \pages 1--37 \endref

\ref \key R2 \bysame
\paper Special positions for surfaces bounded by closed 
braids
\jour Rev. Mat. Iberoamericana
\vol 1 \yr 1985 \pages 93--133 \endref
\endRefs
\enddocument